\documentclass[12pt]{article}
 \usepackage{latexsym}
 \usepackage{amssymb}
 \usepackage{graphicx}
 \usepackage{amsmath}

\newtheorem{Theorem}{Theorem}

\newtheorem{Proposition}{Proposition}

\newtheorem{Corollary}{Corollary}

\newtheorem{Example}{Example}
\newtheorem{Algorithm}{Algorithm}

\newcommand{\A}{{\cal A}}
\newcommand{\B}{{\cal B}}
\newcommand{\C}{{\cal C}}
\newcommand{\D}{{\cal D}}
\newcommand{\I}{{\cal I}}

\newcommand{\uu}{{\bf u}}
\newcommand{\vv}{{\bf v}}
\newcommand{\x}{{\bf x}}
\newcommand{\y}{{\bf y}}

\newcommand{\cc}{{\bf c}}

\newcommand{\0}{{\bf 0}}
\newcommand{\1}{{\bf 1}}
\newcommand{\w}{{\bf w}}

\newcommand{\qed}{\nobreak \ifvmode \relax \else
      \ifdim\lastskip<1.5em \hskip-\lastskip
      \hskip1.5em plus0em minus0.5em \fi \nobreak
      \vrule height0.75em width0.5em depth0.25em\fi}

\def \ep{\hbox{ }\hfill$\Box$}

\begin{document}
\title{Circulant Tensors with Applications to Spectral Hypergraph Theory and Stochastic Process}

\author{Zhongming Chen \thanks{ \footnotesize{School of Mathematical
Sciences and LPMC, Nankai University, Tianjin 300071, P.R. China,
Email: czm183015@mail.nankai.edu.cn.} This author's work was done
when he was visiting The Hong Kong Polytechnic University.} , \quad
Liqun Qi
\thanks{Department of Applied Mathematics, The Hong Kong Polytechnic University, Hung Hom,
Kowloon, Hong Kong. Email: maqilq@polyu.edu.hk. This author's work was supported by the Hong
Kong Research Grant Council (Grant No. PolyU 502510, 502111, 501212
and 501913).}}

\date{\today} \maketitle

\begin{abstract}
\noindent  
Circulant tensors naturally arise from stochastic process and
spectral hypergraph theory.  The joint moments of stochastic
processes are symmetric circulant tensors.  The adjacency, Laplacian
and signless Laplacian tensors of circulant hypergraphs  are also
symmetric circulant tensors.    The adjacency, Laplacian and
signless Laplacian tensors of directed circulant hypergraphs are
circulant tensors, but they are not symmetric in general. In this
paper, we study spectral properties of circulant tensors and their
applications in spectral hypergraph theory and stochastic process.
We show that in certain cases, the largest H-eigenvalue of a
circulant tensor can be explicitly identified. In particular, the
largest H-eigenvalue of a nonnegative circulant tensor can be
explicitly identified.  This confirms the results in circulant
hypergraphs and directed circulant hypergraphs.   We prove that an
even order  circulant B$_0$ tensor is always positive semi-definite.
This shows that the Laplacian tensor and the signless Laplacian
tensor of a directed circulant even-uniform hypergraph are positive
semi-definite. If a stochastic process is $m$th order stationary,
where $m$ is even, then its $m$th order moment, which is a circulant
tensor, must be positive semi-definite. In this paper, we give
various conditions for an even order circulant tensor to be positive
semi-definite.

\noindent {\bf Key words:}\hspace{2mm} Circulant tensors, circulant
hypergraphs, directed circulant hypergraphs, eigenvalues of tensors,
positive semi-definiteness. \vspace{3mm}

\noindent {\bf AMS subject classifications (2010):}\hspace{2mm}
15A18; 15A69
  \vspace{3mm}

\end{abstract}


\section{Introduction}
\hspace{4mm} Circulant matrices are Topelitz matrices.    They form
an important class of matrices in linear algebra and its
applications \cite{Da, HJ, Wik}. As a natural extension of circulant
matrices, circulant tensors naturally arise from stochastic process
and spectral hypergraph theory.  The joint moments of
stochastic processes are symmetric circulant tensors.  The adjacency, Laplacian and signless Laplacian tensors of
circulant hypergraphs  are also symmetric circulant tensors.    The
adjacency, Laplacian and signless Laplacian tensors of directed circulant hypergraphs
are circulant tensors, but they are not symmetric in general.

In this paper, we study spectral properties of circulant tensors and their applications in
spectral hypergraph theory and stochastic process.

Denote $[n] := \{ 1, \cdots, n \}$.  Let $\A = \left(a_{j_1 \cdots
j_m}\right)$ be a real $m$th order $n$-dimensional tensor.    If for
$j_l \in [n-1], l \in [m]$, we have
$$a_{j_1\cdots j_m} = a_{j_1+1\cdots j_m+1},$$
then we say that $\A$ is an $m$th order {\bf Toeplitz tensor}
\cite{BB}. If for $j_l, k_l \in [n], k_l = j_l+1$ mod$(n)$, $l \in
[m]$, we have
\begin{equation} \label{e00}
a_{j_1\cdots j_m} = a_{k_1\cdots k_m},
\end{equation}
then we say that $\A$ is an $m$th order {\bf circulant tensor}.
Clearly, a circulant tensor is a Toeplitz tensor.  By the
definition, all the diagonal entries of a Toeplitz tensor are the
same.   Thus, we may say the {\bf diagonal entry} of a Toeplitz or
circulant tensor.   Tensors with circulant structure were studied in
\cite{RE}.

A real $m$th order $n$-dimensional tensor (hypermatrix) $\A =
(a_{i_1\cdots i_m})$ is a multi-array of real entries $a_{i_1\cdots
i_m}$, where $i_j \in [n]$ for $j \in [m]$. Denote the set of all
real $m$th order $n$-dimensional tensors by $T_{m, n}$.  Then $T_{m,
n}$ is a linear space of dimension $n^m$.   Denote the set of all
real $m$th order $n$-dimensional circulant tensors by $C_{m, n}$.
Then $C_{m, n}$ is a linear subspace of $T_{m, n}$, with dimension
$n^{m-1}$.

Let $\A = (a_{i_1\cdots i_m}) \in T_{m, n}$. If the entries
$a_{i_1\cdots i_m}$ are invariant under any permutation of their
indices, then $\A$ is called a {\bf symmetric tensor}.  Denote the
set of all real $m$th order $n$-dimensional tensors by $S_{m, n}$.
Then $S_{m, n}$ is a linear subspace of $T_{m, n}$.

Let $\A \in T_{m, n}$.  Assume that $m$ is even.  If $\A \x^m \ge 0$
for all $\x \in \Re^n$, then we say that $\A$ is {\bf positive
semi-definite}. If $\A \x^m > 0$ for all $\x \in \Re^n, \x \not =
\0$, then we say that $\A$ is {\bf positive definite}. The
definition of positive semi-definite tensors was first introduced in
\cite{Qi} for symmetric tensors.   Here we extend that definition to
any tensors in $T_{m, n}$.


Throughout this paper, we assume that $m, n \ge 2$.   We use small
letters $x, u, v, \alpha, \cdots$, for scalers, small bold letters
$\x, \y, \uu, \cdots$, for vectors, capital letters $A, B, \cdots$,
for matrices, calligraphic letters $\A, \B, \cdots$, for tensors. We
reserve the letter $i$ for the imaginary unit.    Denote $\1_j \in
\Re^n$ as the $j$th unit vector for $j \in [n]$, $\0$ the zero
vector in $\Re^n$, $\1$ the all $1$ vector in $\Re^n$, and $\hat \1$
the alternative sign vector $(1, -1, 1, -1, \cdots )^\top \in
\Re^n$.  We call a tensor in $T_{m, n}$ the {\bf identity tensor} of $T_{m, n}$, and denote it $\I$
if all of its diagonal entries are $1$ and all of its off-diagonal entries are $0$.

In the next section, we study the applications of circulant tensors
in stochastic process and spectral hypergraph theory.  In
particular, we study what are the concerns of the properties of
circulant tensors in these applications.  If a stochastic process is
$m$th order stationary, where $m$ is even, then its $m$th order
moment, which is a circulant tensor, must be positive semi-definite.
Hence, in the next three sections, we give various conditions for an
even order circulant tensor to be positive semi-definite.

It is well-known that a circulant matrix is generated from the first
row vector of that circulant matrix \cite{Da, HJ, Wik}.    We may
also generate a circulant tensor in this way.    In Section 3, we
define the {\bf root tensor} $\A_1 \in T_{m-1, n}$ and the {\bf
associated tensor} $\bar \A_1 \in T_{m-1, n}$ for a circulant tensor
$\A \in C_{m, n}$.   We show that $\A$ is generated from $\A_1$. It
is also well-known that the eigenvalues and eigenvectors of a
circulant tensor can be written explicitly \cite{Da, HJ, Wik}.   In
Section 3, after reviewing the definitions of eigenvalues and
H-eigenvalues of a tensor in $T_{m, n}$, we show that for any
circulant tensors $\A \in C_{m, n}$ with any $m \ge 2$, including
circulant matrices in $C_{2, n}$, the same $n$ independent vectors
are their eigenvectors. For a circulant tensor $\A \in C_{m, n}$, we
define a one variable polynomial $f_\A(t)$ as its {\bf associated
polynomial}.   Using $f_\A(t)$, we may find the $n$ eigenvalues
$\lambda_k(\A)$ for $k = 0, \cdots, n-1$, corresponding to these $n$
eigenvectors. We call these $n$ eigenvalues the native eigenvalues
of that circulant tensor $\A$. In particular, the first native
eigenvalue $\lambda_0(\A)$, which is equal to the sum of all the
entries of the root tensor, is an H-eigenvalue of $\A$.  We show
that when the associated tensor is a nonnegative tensor,
$\lambda_0(\A)$ is the largest H-eigenvalue of $\A$.  This confirms
the results in circulant hypergraphs and directed circulant
hypergraphs. \textbf{}

In Section 4, we study positive semi-definiteness of an even order
circulant tensor. Recently, it was proved in \cite{QS} that an even
order symmetric B$_0$ tensor is positive semi-definite, and an even
order symmetric B tensor is positive definite.   In Section 4, for
any tensor $\A \in T_{m, n}$, we define a symmetric tensor $\B \in
S_{m, n}$ as its {\bf symmetrization}, and denote it $sym(\A)$.  An
even order tensor is positive semi-definite or positive definite if
and only if its symmetrization is positive semi-definite or positive
definite, respectively.  We show that the symmetrization of a
circulant B$_0$ tensor is still a circulant B$_0$ tensor, and the
symmetrization of a circulant B tensor is still a circulant B
tensor.     This implies that an even order circulant B$_0$ tensor
is always positive semi-definite, and an even order circulant B
tensor is always positive definite. Thus, the Laplacian tensor and
the signless Laplacian tensor of a directed circulant even-uniform
hypergraph are positive semi-definite.   Some other sufficient
conditions for positive semi-definiteness of an even order circulant
tensor are also given in that section.

 In Section 5, we study positive semi-definiteness of even order circulant tensors
 with special root tensors.   When the root tensor $\A_1$ is a diagonal tensor, we show
that in this case, the $n$ native eigenvalues are indeed all the
eigenvalues of that circulant tensor $\A$, with some adequate
multiplicities and more eigenvectors.  We give all such eigenvectors
explicitly.    Then we present some necessary conditions, sufficient
conditions, and necessary and sufficient conditions for an even
order circulant tensor with a diagonal root tensor to be positive
semi-definite.  An algorithm for determining positive
semi-definiteness of an even order circulant tensor with a diagonal
root tensor, and its numerical experiments are also presented. When
the root tensor $\A_1$ itself is a circulant tensor, we call $\A$ a
{\bf doubly circulant} tensor.  We show that when $m$ is even and
$\A_1$ is a doubly circulant tensor itself, if the root tensor of
$\A_1$ is positive semi-definite, then $\A$ is also positive
semi-definite.

\section{Stochastic Process, Circulant Hypergraphs and Directed Circulant Hypergraphs}
\hspace{4mm}
In this section, we study stochastic process, circulant hypergraphs and directed circulant hypergraphs.  We show that circulant tensors naturally arise from these applications.
We study what are the concerns on the properties of circulant tensors in these applications.

\subsection{Stochastic Process}
\hspace{4mm}
For a vector-valued random variable $\x=(x_1, \dots, x_n)$, the joint moment of $\x$ is defined as the expected value of their product:
$$ \text{Mom}(x_{1},\cdots,x_{n})= E\{ x_1 x_2 \cdots x_n \}.  $$
The $m$th order moment of the stochastic vector $\x =(x_1, \dots, x_n)$ is a $m$th order $n$-dimensional tensor, defined by
$$ M_m(\x)= \left[ \text{Mom}(x_{i_1},\cdots,x_{i_m}) \right]_{i_1, \cdots, i_m=1}^n .$$
By definition, we have:
(i) $M_m(\x)$ is symmetric;
(ii) when $m=2$, $M_2(\x)$ is the covariance matrix of the stochastic vector $\x$ with mean $\0$;
(iii) if $\y=A^\top \x$ with $A \in \Re^{n \times N}$, then $ M_m(\y)=M_m(\x) A^m $, where the product is defined in Section \ref{sec3}.

On the other hand, a discrete stochastic process $\x=\{ x_k, k=1,2,\cdots \}$ is called $m$th order stationary if for any points $t_1, \cdots, t_m \in \mathbb{Z}_+$, the joint distribution of $$\{x_{t_1},\cdots, x_{t_m} \}$$ is the same as the joint distribution of $$\{x_{t_1+1},\cdots, x_{t_m+1} \}.$$
A stochastic process is stationary if it is $m$th order stationary for any positive integer $m $. It is well-known that a Markov chain is a stationary process if the initial state is chosen according to the stationary distribution.
We can see that the $m$th order moment of a $m$th order stationary stochastic process $\x$, $M_m(\x)$, is a $m$th order Toeplitz tensor with infinite dimension.
In practice, it may be difficult to handle this case. Instead, a stochastic process $\x=\{ x_k, k=1,2,\cdots \}$ can be approximated by a stochastic process with period $n$, $\x^n = \{ x^n_k, k=1,2,\cdots \}$, where $ x^n_k = x^n_j$ if $k = j$ mod$(n)$. For example,
$\x^1=\{ x_1, x_1, x_1, x_1, \cdots \}$ and $\x^2=\{ x_1, x_2, x_1, x_2, \cdots \}$. We can see that the $m$th order moment of $\x^n$ can be expressed by a $m$th order $n$-dimensional tensor $M_m(\x^n)$ since
$$ \text{Mom}(x^n_{i_1},\cdots,x^n_{i_m})= \text{Mom}(x^n_{j_1},\cdots,x^n_{j_m}),  $$
where $i_k = j_k$ mod$(n)$ for $k \in [m]$.
If the stochastic process $\x$ is $m$th order stationary, the $m$th order moment of the approximation with period $n$, $M_m(\x^n)$, is a circulant tensor of order $m$ and dimension $n$.

\medskip
Given a stochastic process $\x^n$ with period $n$, by Theorem 7.1 of Chapter 9 \cite{Ti}, one can derive that $\x^n$ is the second order stationary if and only if $M_2(\x^n)$ is positive semi-definite. In general, $M_m(\x^n)$ is positive semi-definite when the order $m$ is even.

\begin{Proposition}
For a stochastic process $\x^n$ with period $n$, $M_m(\x^n)$ is positive semi-definite when $m$ is even.
\end{Proposition}
{\bf Proof.} For any $\alpha \in \Re^n$, we have
\begin{eqnarray*}
M_m(\x^n) \alpha^m &=& \sum_{i_1,\cdots,i_m=1}^n \alpha_{i_1} \cdots \alpha_{i_m} \text{Mom}(x^n_{i_1},\cdots,x^n_{i_m})   \\
                   &=& \text{Mom}\left(\sum_{i_1=1}^n \alpha_{i_1} x^n_{i_1},\cdots, \sum_{i_m=1}^n \alpha_{i_m} x^n_{i_m}\right)     \\
                   &=& E\left\{  \left(\sum_{i=1}^n \alpha_i x^n_i\right)^m \right\}.
\end{eqnarray*}
Then, $M_m(\x^n) \alpha^m \geq 0$ since $m$ is even, which means $M_m(\x^n)$ is positive semi-definite.   \ep

This shows that positive semi-definiteness of curculant tensors is important.
In this paper, we will study conditions of positive semi-definiteness of circulant tensors.

\medskip

On the other hand, given a positive semi-definite tensor
$\mathcal{M} \in C_{m,n}$, is there a stationary stochastic process
$\x^n$ with period $n$ such that $M_m(\x^n)=\mathcal{M}$? This
remains as a further studying issue in stochastic process.

\subsection{Circulant Hypergraphs}
\hspace{4mm}    In the recent years, a number of papers appeared in
spectral hypergraph theory via tensors \cite{BP, CD, HQ, HQ1, HQ2,
HQS, LQY, PZ, Qi2, QSW, XC, XC1, XC2}.

A hypergraph $G$ is a pair $(V, E)$, where $V = [n]$ is the set of
vertices and $E$ is a set of subsets of $V$. The elements of $E$ are
called edges. An edge $e \in E$ has the form $e = (j_1, \cdots,
j_m)$, where $j_l \in V$ for $l \in [m]$ and $j_l \not = j_k$ if $l
\not = k$.  The order of $j_1, \cdots, j_m$ is irrelevant for
an edge.  Given an integer $m \geq 2$, a hypergraph $G$ is said to
be $m$-uniform if $|e|=m$ for all $e \in E$, where $|e|$ denotes
number of vertices in the edge $e$.     The degree of a
vertex $j \in V$ is defined as $d(j) = |E(j)|$, where $E(j) = \{ e
\in E: j \in e \}$. If for all $j \in V$, the degrees $d(j)$ have
the same value $d$, then $G$ is called a regular hypergraph, or a
$d$-regular hypergraph to stress its degree $d$.

An $m$-uniform hypergraph $G = (V, E)$ with $V = [n]$ is called a
{\bf circulant hypergraph} if $G$ has the following property: if $e=
(j_1, \cdots, j_m) \in E$, $k_l = j_l +1$ mod$(n), l \in [m]$, then
$\bar e = (k_1, \cdots, k_m) \in E$. Clearly, a circulant hypergraph
is a regular hypergraph.

For an $m$-uniform hypergraph $G = (V, E)$ with $V=[n]$, the
adjacency tensor $\A = \A(G) $ is a tensor in $S_{m, n}$, defined by
$\A = (a_{j_1\cdots j_m})$,
$$
a_{j_1 \cdots j_m}=\frac{1}{(m-1)!}\left\{\begin{array}{ll}
1 &  \text{if } (j_1,\cdots,j_m) \in E \\
0 &  \text{otherwise.}
\end{array}
\right.
$$
The degree tensor $\D=\D(G)$ of $G$, is a diagonal tensor in $S_{m,
n}$, with its $j$th diagonal entry as $d(j)$. The Laplacian tensor
and the signless Laplacian tensor of $G$ are defined by
$\mathcal{L}(G)=\D(G)-\A(G)$ and $\mathcal{Q}(G)=\D(G)+\A(G)$, which
were initially introduced in \cite{Qi2}, and studied further in
\cite{HQ1, HQS, QSW}.  The adjacency tensor, the Laplacian tensor
and the signless Laplacian tensors of a uniform hypergraph are
symmetric.  The adjacency tensor and the signless Laplacian tensor
are nonnegative. The Laplacian tensor and the signless Laplacian
tensor of an even-uniform hypergraph are positive semi-definite
\cite{Qi2}.  It is known \cite{Qi2} that the adjacency tensor, the
Laplacian tensor and the signless Laplacian tensor of a uniform
hypergraph always have H-eigenvalues.    The smallest H-eigenvalue
of the Laplacian tensor is zero with an H-eigenvector $\1$.  The
largest H-eigenvalues of the adjacency tensor and the signless
Laplacian tensor of a $d$-regular hypergraph are $d$ and $2d$
respectively \cite{Qi2}.

Clearly, the adjacency tensor, the Laplacian tensor and the signless
Laplacian tensor of a circulant hypergraph are symmetric circulant
tensors.

\subsection{Directed Circulant Hypergraphs}
\hspace{4mm}   Directed hypergraphs have found applications in
imaging processing \cite{DB}, optical network communications
\cite{LW}, computer science and combinatorial optimization
\cite{GLPN}.   On the other hand, unlike spectral theory of
undirected hypergraphs, it is almost blank for spectral theory of
directed hypergraphs.   In the following, our definition for
directed hypergraphs is the same as in \cite{LW}, which is a special
case of the definition in \cite{DB}, i.e., we discuss the case that
each arc has only one tail.

A directed hypergraph $G$ is a pair $(V, A)$, where $V = [n]$ is the
set of vertices and $A$ is a set of ordered subsets of $V$. The
elements of $A$ are called arcs. An arc $e \in A$ has the form $e =
(j_1, \cdots, j_m)$, where $j_l \in V$ for $l \in [m]$ and $j_l \not
= j_k$ if $l \not = k$.    The order of $j_2, \cdots, j_m$  is
irrelevant.  But the order of $j_1$ is special.   The vertex $j_1$
is called the tail of the arc $e$.   It must be in the first
position of the arc.    The other vertices $j_2, \cdots, j_m$ are
called the heads of the arc $e$.  Similar to $m$-uniform
hypergraphs, we have $m$-uniform directed hypergraphs. The degree of
a vertex $j \in V$ is defined as $d(j) = |A(j)|$, where $A(j) = \{ e
\in A: j \ {\rm is\ a\ tail\ of\ } e \}$.    If for all $j \in V$,
the degrees $d(j)$ have the same value $d$, then $G$ is called a
directed regular hypergraph, or a directed $d$-regular hypergraph.

Similarly, an $m$-uniform directed hypergraph $G = (V, A)$ with $V =
[n]$ is called a {\bf directed circulant hypergraph} if $G$ has the
following property: if $e= (j_1, \cdots, j_m) \in A$, $k_l = j_l +1$
mod$(n), l \in [m]$, then $\bar e = (k_1, \cdots, k_m) \in A$.
Clearly, a directed circulant hypergraph is a regular directed
hypergraph.

For an $m$-uniform directed hypergraph $G = (V, A)$ with $V=[n]$,
the adjacency tensor $\A = \A(G) $ is a tensor in $T_{m, n}$,
defined by $\A = (a_{j_1\cdots j_m})$,
$$
a_{j_1 \cdots j_m}=\frac{1}{(m-1)!}\left\{\begin{array}{ll}
1 &  \text{if } (j_1,\cdots,j_m) \in A \\
0 &  \text{otherwise.}
\end{array}
\right.
$$
Then, the degree tensor $\D=\D(G)$ of $G$, is a diagonal tensor in
$T_{m, n}$, with its $j$th diagonal entry as $d(j)$. The Laplacian
tensor and the signless Laplacian tensor of $G$ are defined by
$\mathcal{L}(G)=\D(G)-\A(G)$ and $\mathcal{Q}(G)=\D(G)+\A(G)$.

The adjacency tensor, the Laplacian tensor and the signless
Laplacian tensors of a uniform directed hypergraph are not symmetric
in general. The adjacency tensor and the signless Laplacian tensor
are still nonnegative. In general, we do not know if the Laplacian
tensor and the signless Laplacian tensor of an even-uniform directed
hypergraph are positive semi-definite or not.  We may still show
that the smallest H-eigenvalue of the Laplacian tensor of an
$m$-uniform directed hypergraph is zero with an H-eigenvector $\1$,
and the largest H-eigenvalues of the adjacency tensor and the
signless Laplacian tensor of a directed $d$-regular hypergraph are
$d$ and $2d$ respectively.

Clearly, the adjacency tensor, the Laplacian tensor and the signless
Laplacian tensor of a directed circulant hypergraph are circulant
tensors.   In general, they are not symmetric.

\section{Eigenvalues of A Circulant Tensor} \label{sec3}
\hspace{4mm} It is well-known that the other row vectors of a
circulant matrix are rotated from the first row vector of that
circulant matrix \cite{Da, HJ, Wik}.    We may also regard a
circulant tensor in this way.   In order to do this, we introduce
{\bf row tensors} for a tensor $\A = (a_{j_1\cdots j_m}) \in T_{m,
n}$.    Let $\A_k = (a^{(k)}_{j_1\cdots j_{m-1}}) \in T_{m-1, n}$ be
defined by $a^{(k)}_{j_1\cdots j_{m-1}} \equiv a_{kj_1\cdots
j_{m-1}}$. We call $\A_k$ the $k$th row tensor of $\A$ for $k \in
[n]$. Let $\A$ be a circulant tensor.   Then we see that the row
tensors $\A_k$ for $k = 2, \cdots, n$, are generated from $\A_1 =
(\alpha_{j_1\cdots j_{m-1}})$, where $\alpha_{j_1\cdots j_{m-1}}
\equiv a^{(1)}_{j_1\cdots j_{m-1}}$. We call $\A_1$ the {\bf root
tensor} of $\A$. We see that $c_0 = \alpha_{1\cdots 1}$ is the
diagonal entry of $\A$. The off-diagonal entries of $\A$ are
generated by the other entries of $\A_1$.  Thus, we define $\bar
\A_1 = (\bar \alpha_{j_1\cdots j_{m-1}}) \in T_{m-1, n}$ by $\bar
\alpha_{1\cdots 1} = 0$ and $\bar \alpha_{j_1\cdots j_{m-1}} =
\alpha_{j_1\cdots j_{m-1}}$ if $(j_1, \cdots, j_{m-1}) \not = (1,
\cdots, 1)$, and call $\bar \A_1$ the {\bf associated tensor} of
$\A$.

We may further quantify this generating operation.   Let $\A =
(a_{j_1\cdots j_m}) \in T_{m, n}$ and $Q=(q_{jk}) \in T_{2, n}$.
Then as in \cite{Qi}, $\B = (b_{k_1\cdots k_m}) \equiv \A Q^m$ is
defined by
$$b_{k_1\cdots k_m} = \sum_{j_1,\cdots, j_m=1}^n a_{j_1\cdots
j_m}q_{j_1k_1}\cdots q_{j_mk_m},$$ for $k_1, \cdots, k_m \in [n]$.
Now we denote $P = (p_{jk}) \in T_{2, n}$ as a permutation matrix
with $p_{jj+1} = 1$ for $j \in [n-1]$, $p_{n1}=1$ and $p_{jk} = 0$
otherwise, i.e.,
\begin{equation} \label{e0}
P = \left(\begin{array}{ccccc} 0 & 1 & \cdots & 0 &
0\\  0 & 0 & 1 & & 0\\ \vdots & 0 & 0 & \ddots & \vdots\\
0 & & \ddots & \ddots & 1\\ 1 & 0 & \cdots & 0 & 0
\end{array} \right).
\end{equation}

 Then, from the definition of circulant tensors, we have the
following proposition.

\begin{Proposition} \label{p1}
Suppose that $\A \in C_{m, n}$ and $P$ is defined by (\ref{e0}).
Then for $k \in [n]$, we have
$$\A_{k+1} = \A_kP^{m-1},$$
where  $\A_{n+1} \equiv \A_1$.
\end{Proposition}

We may also use the definition of circulant tensors to prove the
following proposition.   As the proof is simple, we omit the proof
here.

\begin{Proposition} \label{p2}
Suppose that $\A \in T_{m, n}$ and $P$ is defined by (\ref{e0}).
Then the following three statements are equivalent.

(a). $\A \in C_{m, n}$.

(b). $\A P^m = \A$.

(c). For any $C \in C_{2, n}$, $\A C^m \in C_{m, n}$.
\end{Proposition}

We may denote a circulant matrix $C \in C_{2, n}$ as
\begin{equation} \label{e0.5}
C = \left(\begin{array}{ccccc} c_0 & c_{n-1} & \cdots & c_2 &
c_1\\  c_1 & c_0 & c_{n-1} & & c_2\\ \vdots & c_1 & c_0 & \ddots &
\vdots\\ c_{n-2} & & \ddots & \ddots & c_{n-1}\\ c_{n-1} & c_{n-2} &
\cdots & c_1 & c_0
\end{array} \right).
\end{equation}
It is well-known \cite{Da, HJ, Wik} that the eigenvectors of $\C$
are given by
\begin{equation} \label{e1}
\vv_k = \left(1, \omega_k, \omega_k^2, \cdots,
\omega_k^{n-1}\right)^\top, \end{equation} where $\omega_k = e^{2\pi
ik \over n}$ for $k+1 \in [n]$, with corresponding eigenvalues
$\lambda_k = f_C(\omega_k ),$ where $f_C$ is the associated polynomial
of $C$, defined by
$$f_C(t) = c_0 + c_1t + \cdots + c_{n-1}t^{n-1}.$$
We may also extend this result to circulant tensors.   Note that
$\vv_0 = \1$ is a real vector.

For $\A = (a_{j_1\cdots j_m}) \in T_{m, n}$ and $\x = (x_1, \cdots,
x_n)^\top \in {\boldmath C}^n$, define
$$\A \x^m = \sum_{j_1, \cdots, j_m=1}^n a_{j_1\cdots
j_m}x_{j_1}\cdots x_{j_m}$$  and $\A \x^{m-1}$ as a vector in
${\boldmath C}^n$ with its $j$th component as
$$\left(\A \x^{m-1}\right)_j = \sum_{j_2, \cdots, j_m=1}^n a_{jj_2\cdots
j_m}x_{j_2}\cdots x_{j_m}.$$

Let $\A \in T_{m, n}$.  For any vector $\x \in {\boldmath C}^n$,
$\x^{[m-1]}$ is a vector in ${\boldmath C}^n$, with its $i$th
component as $x_i^{m-1}$.  If $\A \x^{m-1} = \lambda \x^{[m-1]}$ for
some $\lambda \in {\boldmath C}$ and $\x \in {\boldmath C}^n
\setminus \{ \0 \}$, then $\lambda$ is called an eigenvalue of $\A$
and $\x$ is called an eigenvector of $\A$, associated with
$\lambda$. If $\x$ is real, then $\lambda$ is also real.  In this
case, they are called an H-eigenvalue and an H-eigenvector
respectively.   The largest modulus of the eigenvalues of $\A$ is
called the {\bf spectral radius} of $\A$, and denoted as $\rho(\A)$.
The definition of eigenvalues was first given in \cite{Qi} for
symmetric tensors.   It was extended to tensors in $T_{m, n}$ in
\cite{CPZ}.

Suppose that $\A \in C_{m, n}$.  Let its root tensor be $\A_1 =
(\alpha_{j_1\cdots j_{m-1}})$.  Define the {\bf associated
polynomial} $f_\A$ by
\begin{equation} \label{e2}
f_\A(t) = \sum_{j_1, \cdots, j_{m-1}=1}^n \alpha_{j_1\cdots
j_{m-1}}t^{j_1+\cdots+j_{m-1}-m+1}.
\end{equation}

\begin{Theorem} \label{t2}
Suppose that $\A \in C_{m, n}$, its root tensor is $\A_1 =
(\alpha_{j_1\cdots j_{m-1}})$, and its associated tensor is $\bar
\A_1 = (\bar \alpha_{j_1\cdots j_{m-1}})$. Denote the diagonal entry
of $\A$ by $c_0 = a_{1\cdots 1} = \alpha_{1\cdots 1}$.    Then any
eigenvalue $\lambda$ of $\A$ satisfies the following inequality:
\begin{equation} \label{e2.05}
\left|\lambda - c_0\right| \le \sum_{j_1\cdots, j_{m-1}=1}^n
\left|\bar \alpha_{j_1\cdots j_{m-1}}\right|.
\end{equation}
Furthermore, the vectors $\vv_k$, defined by (\ref{e1}), are
eigenvectors of $\A$, with corresponding eigenvalues $\lambda_k =
\lambda_k(\A) = f_\A(\omega_k ),$ where $f_\A$ is the associated
polynomial of $\A$, defined by (\ref{e2}). In particular, $\A$
always has an H-eigenvalue
\begin{equation} \label{e2.1}
\lambda_0 = \lambda_0(\A) = \sum_{j_1,\cdots, j_{m-1}=1}^n
\alpha_{j_1\cdots j_{m-1}},
\end{equation}
with an H-eigenvector $\1$,
and when $n$ is even,
\begin{equation} \label{e2.2}
\lambda_{n \over 2} = \lambda_{n \over 2}(\A) = \sum_{j_1,\cdots,
j_{m-1}=1}^n \alpha_{j_1\cdots j_{m-1}}(-1)^{j_1+\cdots+j_{m-1}-m+1}
\end{equation}
is also an H-eigenvalue of $\A$ with an H-eigenvector $\hat \1$.
\end{Theorem}
\noindent{\bf Proof.}  By the definition of circulant tensors and
Theorem 6(a) of \cite{Qi}, all the eigenvalues of $\A$ satisfy
(\ref{e2.05}). Let $\A_j$ be the $j$th row tensor of $\A$ for $j \in
[n]$. Let $P$ be defined by (\ref{e0}) and $k+1 \in [n]$. It is easy
to verify that $P\vv_k = \omega_k \vv_k$. To prove that $(\vv_k,
\lambda_k)$ is an eigenpair of $\A$, it suffices to prove that for
$j \in [n]$,
\begin{equation} \label{e4}
\A_j \vv_k^{m-1} = \lambda_k \omega_k^{(j-1)(m-1)}. \end{equation}
We prove (\ref{e4}) by induction.   By the definition of the
associate polynomial, we see that (\ref{e4}) holds for $j=1$. Assume
that (\ref{e4}) holds for $j-1$.   By Proposition \ref{p1}, we have
\begin{eqnarray*}
\A_j \vv_k^{m-1} & = & \A_{j-1}P^{m-1}\vv_k^{m-1}\\
& = & \A_{j-1}(P\vv_k)^{m-1} \\
& = & \A_{j-1}(\omega_k\vv_k)^{m-1}\\
& = & \omega_k^{m-1}\A_{j-1} \vv_k^{m-1}\\
& = & \omega_k^{m-1}\lambda_k\omega_k^{(j-2)(m-1)}\\
& = & \lambda_k \omega_k^{(j-1)(m-1)}.
\end{eqnarray*}
This proves (\ref{e4}).  The other conclusions follow from this by the definition of H-eigenvalues and H-eigenvectors.
The proof is completed. \ep

However, unlike a circulant matrix, these $n$ pairs of eigenvalues
and eigenvectors are not the only eigenpairs of a circulant tensor
when $m \ge 3$.  We may see this from the following example.

\begin{Example} \label{ex1}
A circulant tensor $\A = (a_{jkl}) \in C_{3, 3}$ is generated from
the following root tensor
\begin{equation} \label{e5}
\A_1 = \left(\begin{array} {ccc} a & b & c \\ b & c & d \\ c & d & b
\end{array}\right), \end{equation}
where $a= 5.91395$, $b = 2.47255$, $c = 2.92646$, $d = 8.49514$. By
Theorem \ref{t2}, we see that it has eigenvalues $\lambda_0 =
39.1013$, $\lambda_1 = 14.8057+1.1793i$ and $\lambda_2 =
14.8057-1.1793i$.  Using the polynomial system solver {\bf Nsolve}
available in Mathematica, provided by Wolfram Research Inc., Version
8.0, 2010, we may verify that these three eigenvalues are indeed
eigenvalues of $\A$.   However, we found that $\A$ has three more
eigenvalues $\lambda_3 = 4.92535$, $\lambda_4 = -2.08688+13.6795i$
and $\lambda_5 = -2.08688-13.6795i$.
\end{Example}

Thus, for a circulant tensor $\A$, we call the $n$ eigenvalues
$\lambda_k(\A)$ for $k \in [n]$, provided by Theorem \ref{t2}, the
{\bf native eigenvalues} of $\A$, call $\lambda_0(\A)$ the {\bf
first native eigenvalue} of $\A$, and call $\lambda_{n \over 2}(\A)$
the {\bf alternative native eigenvalue} of $\A$ when $n$ is even.

We now show that the first native eigenvalue $\lambda_0(\A)$ plays a
special role in certain cases.

\begin{Theorem} \label{t2.1}
Suppose that $\A \in C_{m, n}$, and its associated tensor is $\bar
\A_1 = (\bar \alpha_{j_1\cdots j_{m-1}})$.   If $\bar \A_1$ is a
nonnegative tensor, then the first native eigenvalue $\lambda_0(\A)$
is the largest H-eigenvalue of $\A$.  If $\bar \A_1$ is a
non-positive tensor, then the first native eigenvalue
$\lambda_0(\A)$ is the smallest H-eigenvalue of $\A$.
\end{Theorem}
\noindent{\bf Proof.}  By Theorem \ref{t2}, we have
$$\lambda_0(\A) = c_0 + \sum_{j_1,\cdots, j_{m-1}=1}^n
\bar \alpha_{j_1\cdots j_{m-1}}.$$ By this and (\ref{e2.05}), the
conclusions hold. \ep

We may apply this theorem to the adjacency, Laplacian and signless
Laplacian tensors of a circulant hypergraph or a directed circulant
hypergraph.  Then we see that the smallest H-eigenvalue of the
Laplacian tensor is zero, the largest H-eigenvalue of the adjacency
tensor is $d$, the largest H-eigenvalue of the signless Laplacian
tensor is $2d$, where $d$ is the common degree of the circulant
hypergraph or the directed circulant hypergraph.   These confirm the
results in Section 2.

When $n$ is even, the alternative native eigenvalue $\lambda_{n
\over 2}(\A)$ also plays a special role in certain cases.  In order
to study the role of the alternative native eigenvalue, we introduce
alternative and negatively alternative tensors. We call a tensor $\B
= (b _{j_1\cdots j_{m}}) \in T_{m, n}$ an {\bf alternative tensor},
if $b_{j_1\cdots j_{m}}(-1)^{\sum_{k=1}^{m} j_k -m} \ge 0$.   We
call $\B$ {\bf negatively alternative} if $-\B$ is alternative.

Then, by definition, we have the following proposition.
\begin{Proposition} \label{p1.1}
Suppose $\B \in T_{m, n} $ and let $\B_k$ be the $k$th row tensor of
$\B$ for $k \in [n]$. Then $\B \in T_{m, n} $ is alternative if and
only if $\B_k$ is alternative when $k$ is odd and $\B_k$ is
negatively alternative when $k$ is even. In particular, $\B_1$ is
alternative if $\B $ is alternative.
\end{Proposition}
\noindent{\bf Proof.} By definition, we have for $k \in [n]$,
$$ b_{k j_1 \cdots j_{m-1}} (-1)^{\sum_{l=1}^{m-1} j_l +k-m} = b_{k j_1 \cdots j_{m-1}} (-1)^{\sum_{l=1}^{m-1} j_l -m+1} (-1)^{k-1} \geq 0 .     $$
It means that when $k$ is odd, $$b^{(k)}_{j_1\cdots j_{m-1}}
(-1)^{\sum_{l=1}^{m-1} j_l -m+1} \geq 0 $$ and when $k$ is even,
$$b^{(k)}_{j_1\cdots j_{m-1}} (-1)^{\sum_{l=1}^{m-1} j_l -m+1} \leq
0 .$$ So the proof is completed.  \ep

However, when $\A$ is circulant, $\A$ may not be alternative even if
$\A_1$ is alternative. A simple counter-example can be given as
follows.

\begin{Example} \label{ex2}
A circulant tensor $\A = (a_{jk}) \in C_{3, 2}$ is given by
\begin{eqnarray*}
\A_1 = \left(\begin{array} {cc} 1 & -1 \\ -1 & 3
\end{array}\right), \quad
\A_2 = \left(\begin{array} {cc} 3 & -1 \\ -1 & 1
\end{array}\right).
\end{eqnarray*}
We can see that $\A_1$ and  $\A_2$ are alternative but by
Proposition \ref{p1.1}, $\A$ is not alternative.
\end{Example}

On the other hand, when $m$ and $n$ are even, we can see that a
circulant tensor is alternative if and only if its root tensor is
alternative.
\begin{Proposition} \label{p1.2}
Suppose $\A \in C_{m,n}$, where $m$ and $n$ are even. Then, $\A$ is
(negatively) alternative if and only if its root tensor $\A_1$ is
(negatively) alternative.
\end{Proposition}
\noindent{\bf Proof.} By Proposition \ref{p1.1}, we only prove that
$\A$ is alternative if its root tensor $\A_1$ is alternative. Let
$\A_k$ be the $k$th row tensor of $\A$ for $k \in [n]$. We first
show that $\A_2$ is negatively alternative since $\A_1$ is
alternative. For any $j_1, \cdots, j_{m-1} \in [n]$, let $s$ be the
number of the indexes that are equal to $1$. Without loss of
generality, we assume $j_1=\cdots=j_s=1$. By Proposition \ref{p1},
we have
\begin{eqnarray*}
a^{(2)}_{j_1\cdots j_{m-1}} (-1)^{\sum_{l=1}^{m-1} j_l -m+1} & = & a^{(2)}_{j_1\cdots j_{m-1}} (-1)^{\sum_{l=s+1}^{m-1}(j_l-1)}   \\
                                                             & = & a^{(1)}_{n \cdots n j_{s+1}-1 \cdots j_{m-1}-1} (-1)^{\sum_{l=s+1}^{m-1}(j_l-1)} \\
                                                             & = & a^{(1)}_{n \cdots n j_{s+1}-1 \cdots j_{m-1}-1} (-1)^{ns+\sum_{l=s+1}^{m-1}(j_l-1)-m+1} (-1)^{m-1-ns}    \\
                                                             & \leq & 0 .
\end{eqnarray*}
The last inequality holds because $\A_1$ is alternative and $m-1-ns$
is odd for any $s \in [n]\cup \{ 0 \}$ since $m$ and $n$ are even.
By induction, one can obtain that $\A_k$ is alternative when $k$ is
odd and $\A_k$ is negatively alternative when $k$ is even, which
means that $\A$ is alternative by Proposition \ref{p1.1}. \ep

\begin{Theorem} \label{t2.2}
Let $n$ be even.   Suppose that $\A \in C_{m, n}$, and its
associated tensor is $\bar \A_1 = (\bar \alpha_{j_1\cdots
j_{m-1}})$.   If $\bar \A_1$ is an alternative tensor, then the
alternative native eigenvalue $\lambda_{n \over 2}(\A)$ is the
largest H-eigenvalue of $\A$.  If $\bar \A_1$ is a negatively
alternative tensor, then the alternative native eigenvalue
$\lambda_{n \over 2}(\A)$ is the smallest H-eigenvalue of $\A$.
\end{Theorem}
\noindent{\bf Proof.}  Let $n$ be even.   By Theorem \ref{t2}, we
have
$$\lambda_{n \over 2}(\A) = c_0 + \sum_{j_1\cdots, j_{m-1}=1}^n
\bar \alpha_{j_1\cdots j_{m-1}}(-1)^{\sum_{k=1}^{m-1} j_k-m+1}.$$ By
this and (\ref{e2.05}), the conclusions hold. \ep

\bigskip

Note that the native eigenvalues other than $\lambda_0(\A)$ and $\lambda_{n \over 2}(\A)$ are in general not H-eigenvalues.

\section{Positive Semi-definiteness of Even Order Circulant Tensors}\label{sec4}
\hspace{4mm}   Let $j_l \in [n]$ for $l \in [m]$.  Define the
generalized Kronecker symbol \cite{Qi, QS} by
$$
\delta_{j_1 \cdots j_m}=\left\{\begin{array}{ll}
1 &  \text{if } j_1=\cdots = j_m \\
0 &  \text{otherwise.}
\end{array}
\right.
$$

Suppose that $\A = (a_{j_1\cdots j_m}) \in T_{m, n}$.   We say that
$\A$ is a B$_0$ tensor if for all $j \in [n]$
\begin{equation} \label{b1}
\sum_{j_2,\cdots, j_m=1}^n a_{jj_2\cdots j_m} \ge 0
\end{equation} and
\begin{equation} \label{b2}
{1 \over n^{m-1}}\sum_{j_2,\cdots, j_m=1}^n a_{jj_2\cdots j_m} \ge
a_{jk_2\cdots k_m}, \ {\rm if\ } \delta_{jk_2\cdots k_m} = 0.
\end{equation}
If strict inequalities hold in (\ref{b1}) and (\ref{b2}), then $\A$
is called a B tensor \cite{QS}.  The definitions of B and B$_0$ tensors are generalizations of the definition of B matrix \cite{Pe}. It was proved in \cite{QS} that
an even order symmetric B tensor is positive definite and an even
order symmetric B$_0$ tensor is positive semi-definite.   We may
apply this result to even order symmetric circulant B$_0$ or B
tensors. What we wish to show is that an even order circulant B
tensor is positive definite and an even order circulant B$_0$ tensor
is positive semi-definite, i.e., we do not require the tensor to be
symmetric here.   In this way, we may apply our result to directed
circulant hypergraphs.    The tool for realizing this is
symmetrization.

By the definition of circulant tensors, it is easy to see that for
$\A = (a_{j_1\cdots j_m}) \in C_{m, n}$, $\A$ is a circulant B$_0$
tensor if and only if
\begin{equation} \label{b3}
\sum_{j_1,\cdots, j_m=1}^n a_{j_1\cdots j_m} \ge 0
\end{equation} and
\begin{equation} \label{b4}
{1 \over n^m}\sum_{j_1,\cdots, j_m=1}^n a_{j_1\cdots j_m} \ge \max
\{ a_{k_1\cdots k_m} : \delta_{k_1\cdots k_m} = 0 \}.
\end{equation}
If strict inequalities hold in (\ref{b3}) and (\ref{b4}), then $\A$
is a circulant B tensor.

It was established in \cite{Qi} that an even order real symmetric
tensor has always H-eigenvalues, and it is positive semi-definite
(positive definite) if and only if all of its H-eigenvalues are
nonnegative (positive). This is not true in general for a
non-symmetric tensor.   In order to use the first native eigenvalue
or the alternative eigenvalue of a nonsymmetric circulant tensor to
check its positive semi-definiteness, we may also use
symmetrization.

We now link a general tensor $\A \in T_{m, n}$ to a symmetric tensor
$\B \in S_{m, n}$.

Let $\A \in T_{m, n}$.  Then there is a unique symmetric tensor $\B
\in S_{m, n}$ such that for all $\x \in \Re^n$, $\A \x^m = \B \x^m$.
We call $\B$ the {\bf symmetrization} of $\A$, and denote it
$sym(\A)$.  Thus, when $m$ is even, a tensor $\A \in T_{m, n}$ is
positive semi-definite (positive definite) if and only if all of the
H-eigenvalues of $sym(\A)$ are nonnegative (positive).

We call an index set $(k_1, \cdots, k_m)$ a permutation of another index set ${j_1, \cdots, j_m}$ if $\{ k_1, \cdots, k_m \} = \{ j_1, \cdots, j_m \}$, denote this operation by $\sigma$,
and denote $\sigma(j_1,\cdots, j_m) = (k_1, \cdots, k_m)$.   Denote the set of all distinct permutations of an index set $(j_1, \cdots, j_m)$, by $\Sigma(j_1, \cdots, j_m)$.  Note that $|\Sigma(j_1, \cdots, j_m)|$, the cardinality of $\Sigma(j_1, \cdots, j_m)$, is variant for different index sets.   For example, if $j_1=\cdots = j_m$, then $|\Sigma(j_1, \cdots, j_m)|=1$; but if
all of $j_1, \cdots, j_m$ are distinct, $|\Sigma(j_1, \cdots, j_m)|=m!$.

Let $\A = (a_{j_1\cdots j_m}) \in T_{m, n}$ and $sym(\A) = \B = (b_{j_1\cdots j_m})$.   Then it is not difficult to see that
\begin{equation} \label{e3.1}
b_{j_1\cdots j_m} = { \sum_{\sigma \in \Sigma(j_1, \cdots, j_m)} a_{\sigma(j_1, \cdots, j_m)} \over |\Sigma(j_1, \cdots, j_m)|}.
\end{equation}

For any $A \in T_{m, n}$, we use $\D(\A)$ to denote a diagonal
tensor in $T_{m, n}$, whose diagonal entries are the same as those
of $\A$.

With this preparation, we are now ready to prove the following theorem.

\begin{Theorem} \label{t3.1}
Let $\A = (a_{j_1\cdots j_m}) \in T_{m, n}$.   Then  we have the
following conclusions:

(a). $\D(\A) = \D(sym(\A))$.

(b). If $\A -\D(\A)$ are nonnegative (or non-positive or alternative
or negatively alternative, respectively), then $sym(\A)
-\D(sym(\A))$ are also nonnegative (or non-positive or alternative
or negatively alternative, respectively).

(c). The symmetrization of a Toeplitz tensor is still a Toeplitz
tensor. The symmetrization of a circulant tensor is still a
circulant tensor.

(d).  The symmetrization of a circulant B$_0$ tensor is still a
circulant B$_0$ tensor.  The symmetrization of a circulant B tensor
is still a circulant B tensor.

(e). Suppose that $\A \in C_{m, n}$.  Then we have
$$\lambda_0(\A) = \lambda_0(sym(\A)).$$
If the associated tensor of a circulant tensor is
 nonnegative (or non-positive), then the associated tensor of the symmetrization of a circulant tensor is
 also nonnegative (or non-positive).

(f). Suppose $\A \in C_{m,n}$, where $m$ and $n$ are even. Then, we
have
$$\lambda_{\frac{n}{2}} (\A)  = \lambda_{\frac{n}{2}} ( sym(\A) ).$$
 \end{Theorem}

\noindent{\bf Proof.}  We have (a) and (b) from (\ref{e3.1})
directly.

(c). Let $\A = (a_{j_1\cdots j_m}) \in T_{m, n}$ be a Toeplitz
tensor, and  $sym(\A) = \B = (b_{j_1\cdots j_m})$. By (\ref{e3.1}),
for $j_l \in [n-1], l \in [m]$,
\begin{eqnarray*}
b_{j_1\cdots j_m} & = & { \sum_{\sigma \in \Sigma(j_1, \cdots, j_m)} a_{\sigma(j_1, \cdots, j_m)} \over |\Sigma(j_1, \cdots, j_m)|} \\
& = & { \sum_{\sigma \in \Sigma(j_1+1, \cdots, j_m+1)} a_{\sigma(j_1+1, \cdots, j_m+1)} \over |\Sigma(j_1+1, \cdots, j_m+1)|} \\
& = & b_{j_1+1\cdots j_m+1}.
\end{eqnarray*}
Thus, $sym(\A) = \B$ is a Toeplitz tensor.   When $\A$ is a
circulant tensor, we may prove that $sym(\A)$ is a circulant tensor
similarly.

(d).  Let $\A = (a_{j_1\cdots j_m}) \in C_{m, n}$ and $sym(\A) = \B
= (b_{j_1\cdots j_m})$.   By (c), $\B \in C_{m, n}$.   Suppose now
that $\A$ ia B$_0$ tensor.   By (\ref{e3.1}) and (\ref{b3}), we have
$$\sum_{j_1,\cdots, j_m=1}^n b_{j_1\cdots j_m} = \sum_{j_1,\cdots, j_m=1}^n a_{j_1\cdots j_m} \ge
0.$$ By (\ref{e3.1}) and (\ref{b4}), we have
$${1 \over
n^m}\sum_{j_1,\cdots, j_m=1}^n b_{j_1\cdots j_m} \ge {1 \over
n^m}\sum_{j_1,\cdots, j_m=1}^n a_{j_1\cdots j_m} \ge \max \{
a_{k_1\cdots k_m} : \delta_{k_1\cdots k_m} = 0 \} \ge \max \{
b_{k_1\cdots k_m} : \delta_{k_1\cdots k_m} = 0 \}.$$ Thus, $\B$ is
also a B$_0$ tensor.   Similarly, If $\A$ is a B tensor, then $\B$
is also a B tensor.

(e).  Let $\A \in C_{m, n}$.   The equality $\lambda_0(\A)  =
\lambda_0(sym(\A))$ holds because
$$
\lambda_0(\A)  = \frac{1}{n} \A \vv_0^{m} = \frac{1}{n} sym(\A)
\vv_0^{m} = \lambda_0(sym(\A)).
$$
The last conclusion follows from (b).

(f). For $ k+1 \in [n]$, let $\omega_k$ and $\vv_k$ be defined in
(\ref{e1}). By Theorem \ref{t2}, one can obtain
$$
\lambda_k(\A) \vv_k^\top \vv_k^{[m-1]}= \A \vv_k^m = sym(\A) \vv_k^m = \lambda_k(sym(\A)) \vv_k^\top \vv_k^{[m-1]}.
$$
By simple computation, we have
$$
\vv_k^\top \vv_k^{[m-1]}= \sum_{j=1}^n \omega_k^{m(j-1)} = \left\{ \begin{array}{ll}
                                                                       \frac{1-\omega_k^{nm}}{1-\omega_k^{m}} =0 &  \text{if} \ \omega_k^{m} \neq 1 ,  \\
                                                                        m                                        &  \text{if} \ \omega_k^{m} =
                                                                        1.
                                                                       \end{array}
                                                                \right.
$$
In particular, since $m$ and $n$ are even, we have $\omega_{\frac{n}{2}}^m = (-1)^m=1 $. It follows that $\lambda_{\frac{n}{2}} (\A)  = \lambda_{\frac{n}{2}} ( sym(\A) )$. \ep

\medskip

In fact, from the proof of Theorem \ref{t3.1}, we can see that $\lambda_k (\A)  = \lambda_k ( sym(\A) )$ if $\omega_k^{m} = 1$. And the equality $\lambda_0(\A) = \lambda_0(sym(\A))$ holds since $\omega_0 = 1 $.   Note that when $m$ is odd, the equality $\lambda_{\frac{n}{2}} (\A)  =
\lambda_{\frac{n}{2}} ( sym(\A) )$  may not hold. See Example \ref{ex2}. By computation, $ sym(\A)
\in C_{3,2} $ is generated by the root tensor
$$
sym(\A)_1 = \left(\begin{array} {cc} 1 & 1/3 \\ 1/3 & 1/3
\end{array}\right).
$$
We can see that $ \lambda_1 (\A) = 6 $ and $ \lambda_1 (sym(\A)) = 2/3 $.
On the other hand, we can also see that $ \lambda_0 (\A) =  \lambda_0 (sym(\A)) = 2 $
and $ \lambda_1 (\A) $ is the largest H-eigenvalue
of $ \A $ since $\A_1$ is alternative.

\medskip

We now have the following corollaries.

\begin{Corollary} \label{ce-1} An even order circulant B$_0$ tensor
is positive semi-definite.    An even order circulant B tensor is
positive definite.
\end{Corollary}
\noindent{\bf Proof.}  Suppose that $\A$ is an even order circulant
B$_0$ tensor.  Then by (d) of Theorem \ref{t3.1}, $\B =$ sym$(\A)$
is also an even order circulant B$_0$ tensor.   Since $\B$ is
symmetric, by \cite{QS}, it is positive semi-definite.   Since $\A$
is positive semi-definite if and only if sym$(\A)$ is positive
semi-definite.    The other conclusion holds similarly. \ep

Note that an even order B$_0$ tensor may not be positive
semi-definite. Let
$$A =\left(\begin{array}{cc} 10 & 10\\  1 & 1
\end{array} \right).$$  Then $A$ is a B$_0$ tensor.    Let $\x = (1, -9)^\top$.   Then $\x^\top A
\x = - 8$. Thus, $A$ is not positive semi-definite.

In the next corollary, we stress that we may use (\ref{b3}) and
(\ref{b4}) instead of (\ref{b1}) and (\ref{b2}) to check an even
order circulant tensor is positive semi-definite or not.   The
conditions (\ref{b3}) and (\ref{b4}) contain less number of
inequalities than (\ref{b1}) and (\ref{b2}).

\begin{Corollary}  \label{ce-2}   Suppose that $\A = (a_{j_1\cdots j_m}) \in C_{m,
n}$ and $m$ is even.  If (\ref{b3}) and (\ref{b4}) hold, then $\A$
is positive semi-definite.  If strict inequalities hold in
(\ref{b3}) and (\ref{b4}), then $\A$ is positive definite.
\end{Corollary}

We may apply these two corollaries to directed circulant
hypergraphs.

\begin{Corollary}  \label{ce-3}   The Laplacian tensor and the
signless Laplacian tensor of a directed circulant even-uniform
hypergraph are positive semi-definite.
\end{Corollary}

As positive semi-definiteness of the Laplacian tensor and the
signless Laplacian tensor of an even-uniform hypergraph plays an
important role in spectral hypergraph theory \cite{HQ, HQ1, HQ2,
HQS, LQY, Qi2, XC, XC2}, The above result will be useful in the
further research for directed circulant hypergraphs.

We may have some other corollaries of Theorem \ref{t3.1} as follows.

\begin{Corollary}  \label{ce1} Suppose that $m$ is even.  If the associated
tensor of a circulant tensor $\A$ is non-positive, then $\A$ is
positive semi-definite if and only if $\lambda_0(\A)$ is
nonnegative.
\end{Corollary}

\begin{Corollary} \label{ce2}
Suppose that $m$ and $n$ are even, $\A \in C_{m, n}$, and its associate tensor $\bar \A_1$
is negatively alternative.    Then $\A$ is positive semi-definite if
and only if $\lambda_{n \over 2}(\A) \ge 0$.
\end{Corollary}
\noindent{\bf Proof.} By definition, we can see that the associate tensor $\bar \A_1$ is the root tensor of $\A - \D(\A)$. By Proposition \ref{p1.2}, one can derive that $\A - \D(\A)$ is negatively alternative since $m$ and $n$ are even. By Theorem \ref{t3.1} (b), it follows that $sym(\A) - \D(sym(\A))$ is also
negatively alternative. Again, by Proposition \ref{p1.2}, $\bar {sym(\A)}_1$ is also negatively alternative.  By Theorem \ref{t2.2}, in this case, $\lambda_{n
\over 2}(sym(\A))$ is the smallest H-eigenvalue of $sym(\A)$. By Theorem \ref{t3.1} (e),
we have $\lambda_{\frac{n}{2}} (\A)  = \lambda_{\frac{n}{2}} ( sym(\A) )$.  The
conclusion follows now. \ep

\begin{Corollary}  \label{ce3} Suppose that $m$ is even.  Suppose that $\A \in C_{m, n}$ is positive semi-definite, and its diagonal entry is $c_0$.
Then $c_0 \ge 0$ and $\lambda_0(\A) \ge 0$.    If furthermore that
$n$ is even, then  $\lambda_{n \over 2}(\A) \ge 0$.
\end{Corollary}


\bigskip

We may further establish a sufficient condition for positive
semi-definiteness of an even order circulant tensor.   For any
tensor $\A = (a_{j_1\cdots j_m}) \in T_{m, n}$, we denote $|\A|
\equiv (|a_{j_1\cdots j_m}|) \in T_{m, n}$.

\medskip

\begin{Theorem} \label{t7.0}
Suppose that $m$ ie even, $\A = (a_{j_1\cdots j_m}) \in C_{m, n}$,
the diagonal entry of $\A$ is $c_0$, and the associated tensor of
$\A$ is $\bar \A_1 = (\bar \alpha_{j_1\cdots j_{m-1}})$.  If
\begin{equation} \label{e3.2}
c_0 \ge \sum_{j_1, \cdots, j_{m-1}=1}^n |\alpha_{j_1\cdots
j_{m-1}}|,
\end{equation}
then $\A$ is positive semi-definite.
\end{Theorem}
\noindent{\bf Proof.} Suppose that (\ref{e3.2}) holds.   For any $\x
\in \Re^n$,
\begin{eqnarray*}
\A \x^m & =  &\D(\A)\x^m + (\A - \D(\A))\x^m \\
& \ge & c_0 \|\x\|_m^m - |\A - \D(\A)| \x^m \\
& \ge & c_0 \|\x\|_m^m - \lambda_0(|\A - \D(\A)|)\|\x \|_m^m \\
& = & c_0 \|\x\|_m^m - \left[\sum_{j_1, \cdots, j_{m-1}=1}^n |\alpha_{j_1\cdots j_{m-1}}|\right] \|\x \|_m^m \\
& \ge & 0.
\end{eqnarray*}
Here, the second inequality holds because that $|\A - \D(\A)|$ is a
nonnegative circulant tensor, hence its spectral radius is equal to
its largest H-eigenvalue  $\lambda_0(|\A - \D(\A)|)$ by Theorem
\ref{t2.1}.   It is easy to see that $|\bar \A_1|$ is the root
tensor of $|\A - \D(\A)|$.   Then we have the second equality.  The
third inequality follows from (\ref{e3.2}).    This shows that $\A$
is positive semi-definite. \ep

Note that Corollary \ref{ce-2} does not imply Theorem \ref{t7.0},
and Theorem \ref{t7.0} does not imply Corollary \ref{ce-2}.  Thus,
they are two different sufficient conditions for positive
semi-definiteness of even order circulant tensors.

\section{Circulant Tensors with Special Root Tensors}
\hspace{4mm}    In this section, we consider conditions for positive
semi-definiteness of even order circulant tensors with special root
tensors, including diagonal root tensors and circulant root tensors.

\subsection{Circulant Tensors with Diagonal Root Tensors}
\hspace{4mm}  Suppose that $\A \in C_{m, n}$ and $\A_1$ is its root tensor.   Assume that $\A_1 = (\alpha_{j_1\cdots j_{m-1}})$ is a diagonal tensor, with
$\alpha_{j_1\cdots j_{m-1}} = c_{j-1}$ if $j_1 = \cdots = j_{m-1} = j \in [n]$, and $\alpha_{j_1\cdots j_{m-1}} = 0$ otherwise.
In this case, we may give all the eigenvalues and eigenvectors (up to some scaling constants) explicitly.  Such a circulant tensor may be one of the simple cases of circulant tensors.   We study its properties such that we can understand more about circulant tensors.

\begin{Theorem} \label{t3}
Let circulant matrix $C$ be defined by (\ref{e0.5}).  With the above
assumptions, the $n$ native eigenvalues $\lambda_k$ of $\A$ are all
possible eigenvalues of $\A$.   They are exactly the $n$ eigenvalues
of the circulant matrix $C$.  For $k+1 \in [n]$,  each eigenvalue
$\lambda_k$ has the following eigenvectors $\y_{kl} = (1, \eta_{kl},
\eta_{kl}^2, \cdots, \eta_{kl}^{n-1})^\top$, where $\eta_{kl} =
e^{2\pi kli \over n(m-1)}$ for $l+1 \in [m-1]$.
\end{Theorem}
\noindent{\bf Proof.}   Let $\y = (y_1, \cdots, y_n)^\top \in {\boldmath C}^n
\setminus \{ \0 \}$ and $\lambda$ be an eigenpair of $\A$.  Define $c_{j-n} = c_j$ for $j \in [n]$.  Then for $j \in [n]$, we have \begin{equation} \label{e7}
\lambda y_j^{m-1} = \left(\A \y^{m-1}\right)_j = \sum_{l=1}^n c_{l-j}y_l^{m-1}.\end{equation}
Let $\x = (y_1^{m-1}, y_2^{m-1}, \cdots, y_n^{m-1})^\top$.   Then we see that (\ref{e7}) is equivalent to $\lambda \x = C \x$, i.e., $(\lambda, \x)$ form an eigenpair of circulant matrix $C$. Now the conclusion can be derived easily. \ep

It is easy to see that $\A$, the circulant tensor with a diagonal root tensor discussed above, is symmetric if and only if $c_j = 0$ for $j \in [n-1]$.    Thus, in general, such a circulant tensor is not symmetric.

Now we discuss positive semi-definiteness of a circulant tensor with a diagonal root tensor.  First, by direct derivation, we have the following result.

\begin{Proposition} \label{p3}
Let $\A \in C_{m,n}$ have a diagonal root tensor as described above.   Then for any $\x = (x_1, \cdots, x_n)^\top \in \Re^n$,
\begin{equation} \label{e8}
\A \x^m = \x^\top C^\top \x^{[m-1]} \equiv \sum_{j, l = 1}^n
c_{l-j}x_jx_l^{m-1} = c_0\sum_{l=1}^n x_l^m + \sum_{j, l = 1 \atop j
\not = l}^n c_{l-j}x_jx_l^{m-1},
\end{equation}
where $C$ is the circulant matrix defined by (\ref{e0.5}).
\end{Proposition}

\begin{Example} \label{ex3}
Let $m =4$ and $n=2$.   Let $c_0 = c_1 = 1$.  Then by Proposition \ref{p3},$$\A \x^4 = x_1^4 + x_1^3x_2 + x_1x_2^3 + x_2^4 = (x_1+x_2)^2(x_1^2 -x_1x_2 + x_2^2) \ge 0$$ for any $\x \in \Re^2$.  Thus, $\A$ is positive semi-definite.
\end{Example}

By (\ref{e2.1}) and (\ref{e2.2}), we have
\begin{equation} \label{e8.1}
\lambda_0(\A) = \sum_{j=0}^{n-1} c_j,
\end{equation}
and when $n$ is even,
$$\lambda_{n \over 2}(\A) = \sum_{j=0}^{n-1} c_j(-1)^{j(m-1)}.$$
In particular, when $m$ is also even, we have
\begin{equation} \label{e8.2}
\lambda_{n \over 2}(\A) = \sum_{j=0}^{n-1} c_j(-1)^j.
\end{equation}

Let $\cc = (c_1, \cdots, c_{n-1})^\top \in \Re^{n-1}$.  Let $k \le {n \over 2}$.    We say that $\cc$ is {\bf $k$-alternative} if $n = 2pk$ for some integer $p$, $c_{(2q-1)k} \ge 0$ and $c_{2qk} \le 0$ for $q \in [p]$ and $c_j = 0$ otherwise.   When $n = 2pk$, let
$\hat \1^{(k)}$ be a vector in $\Re^n$ such that $\hat \1^{(k)}_j = 1$ for $(2q-2)k+1 \le j \le (2q-1)k$ and $\hat \1^{(k)}_j = -1$ for $(2q-1)k+1 \le j \le 2qk$, for $q \in [p]$.

In the following, we give some necessary conditions, sufficient conditions, necessary and sufficient conditions for an even order circulant tensor with a diagonal root tensor to be positive semi-definite.

\begin{Theorem} \label{t6}
Let $\A \in C_{m,n}$ have a diagonal root tensor as described at the beginning of this section.
Suppose that $m$ is even.    Then, we have the following conclusions:

(a). If $\A$ is positive semi-definite, then $c_0 \ge 0$ and
$\lambda_0(\A) \ge 0$.  If furthermore $n$ is even, then $\lambda_{n
\over 2}(\A) \ge 0$.

(b).  If
\begin{equation} \label{e9}
c_0 \ge \sum_{j=1}^{n-1} |c_j|,
\end{equation}
then $\A$ is positive semi-definite.

(c).  If $\cc$ is non-positive, then $\A$ is positive semi-definite if and only if (\ref{e9}) holds.

(d).  If $n = 2pk$ for some positive integers $p$ and $k$, and $\cc$ is $k$-alternative, then $\A$ is positive semi-definite if and only if (\ref{e9})
holds.
\end{Theorem}
\noindent{\bf Proof.} (a).   This follows from Corollary \ref{ce3}.

(b). This follows from Theorem \ref{t7.0}.

(c). If $\cc$ is non-positive, then the associated tensor of $\A$ is non-positive.    By Corollary \ref{ce1}, $\A$ is positive semi-definite if and only if
$$\lambda_0(\A) = \sum_{j=0}^n c_j \ge 0.$$
Since $\cc$ is non-positive and $c_0 \ge 0$, the above inequality holds if and only if (\ref{e9}) holds.    This proves (c).

(d). Suppose that $n = 2pk$ for some positive integers $p$ and $k$, and $\cc$ is $k$-alternative.   Then (\ref{e9}) holds in this case.    By (b), $\A$ is positive semi-definite.   On the other hand, if (\ref{e9}) does not hold, Let $\x = \hat \1^{(k)}$ in (\ref{e8}).  We have $\A \x^m < 0$, i.e., $\A$ is not positive semi-definite.  This proves (d).  The theorem is proved. \ep

\bigskip

Are there some other cases such that (\ref{e9}) is also a sufficient and necessary condition such that $\A$ is positive semi-definite?

\medskip

Suppose that $m$ is even. Can we give all the H-eigenvalues of $sym(\A)$
explicitly?  If so, we may determine $\A$ is positive semi-definite
or not.   Otherwise, can we construct an algorithm to find
the global optimal value of one of the following two minimization
problems when $m$ is even? The two minimization problems are as
follows:
\begin{eqnarray} \label{e8.8}
& \min & c_0\sum_{l=1}^n x_l^m + \sum_{j, l = 1 \atop j \not = l}^n
c_{l-j}x_jx_l^{m-1} \nonumber \\
& {\rm subject\ to}\ & \sum_{j=1}^nx_j^2 = 1,
\end{eqnarray}
and
\begin{eqnarray} \label{e8.9}
& \min & c_0\sum_{l=1}^n x_l^m + \sum_{j, l = 1 \atop j \not = l}^n
c_{l-j}x_jx_l^{m-1} \nonumber \\
& {\rm subject\ to}\ & \sum_{j=1}^nx_j^m = 1.
\end{eqnarray}
By Proposition \ref{p3}, $\A$ is positive semi-definite if and only
if the global optimal value of (\ref{e8.8}) or (\ref{e8.9}) is
nonnegative.   In Subsection 5.3, we will give an algorithm to determine positive semi-definiteness of an even order circulant tensor with a diagonal root tensor.

\subsection{Doubly Circulant Tensors}
\hspace{4mm} Let $\A \in C_{m, n}$.   If its root tensor $\A_1$
itself is a circulant tensor, by Propositions \ref{p1} and \ref{p2},
we see that all the row tensors of $\A$ are duplicates of $\A_1$,
i.e., $\A_k = \A_1$ for $k \in [n]$.    We call such a circulant
tensor $\A$ a {\bf doubly circulant tensor}.

Let $\A$ be an even order doubly circulant tensor.  Suppose that
$\A_{11} \in T_{m-2,n}$ is the root tensor of $\A_1$.  A natural question is that if there is a
relation between $\A_{11}$ and $\A$ in terms of the positive
semi-definiteness, i.e., if $\A_{11}$ is positive semi-definite, is
$\A$ also positive semi-definite?  And if $\A$ is positive
semi-definite, is $\A_{11}$ also positive semi-definite?
Unfortunately, the answers to these two questions are both ``no''.
See the following example.

\begin{Example} \label{ex4}
Let $\A_{11} = diag\{d_1, d_2 \} $. Then, for any $\x \in \Re^2 $,
we have
\begin{eqnarray*}
\A \x^4 &= & (x_1+x_2) \A_1 \x^3   \\
        &= & (x_1+x_2) [d_1 (x_1^3 +x_2^3) + d_2 (x_1^2 x_2 + x_2^2 x_1)] \\
        &= & (x_1+x_2)^2 [ d_1 x_1^2 +(d_2-d_1) x_1 x_2+ d_1 x_2^2 ].
\end{eqnarray*}

\noindent \textbf{Case 1:} $d_1=1$, $d_2=5$. $\A_{11}$ is positive
semi-definite. However, $\A$ is not positive semi-definite since $\A
\x^4 <0$ for $\x=(1,  -2)^\top$.

\noindent \textbf{Case 2:} $d_1=1$, $d_2=-0.5$. $\A$ is positive
semi-definite since $\A \x^4  =(x_1+x_2)^2  [ x_1^2 - 1.5 x_1 x_2+
x_2^2 ] \geq 0$. However, $\A_{11}$ is not positive semi-definite
since $d_2<0$.
\end{Example}

However, we may answer this question positively if $\A_1$ is also doubly circulant.
\begin{Proposition}
If $\A \in T_{m,n}$ is a doubly circulant tensor and $\A_1$ is its root tensor, then for any $\x \in \Re^n$, we have
$$ \A \x^m = \sum_{k=1}^n x_k \A_1 \x^{m-1}. $$
On the other hand, if $m$ is even and $\A_1$ is doubly circulant, then we have the following conclusions:

(a). $\A$ is doubly circulant.

(b). If $\A_{11}$ is positive semi-definite, then $\A$ is positive semi-definite.

(c). If $\A$ is positive semi-definite, then for any $\x \in \Re^n$ satisfying $\sum_{k=1}^n x_k \neq 0 $, we have
$$\A_{11} \x^{m-2} \geq 0 ,$$
where $\A_{11}$ is the root tensor of $\A_1$.
\end{Proposition}
\noindent{\bf Proof.} If $\A$ is doubly circulant, then we have $\A_k = \A_1$ for $k \in [n]$. It follows that for any $\x \in \Re^n$, one can obtain
$$ \A \x^m = \sum_{k=1}^n x_k \A_1 \x^{m-1}.          $$
On the other hand, if $\A_1$ is doubly circulant, then we have $\A$ is doubly circulant by definition and
$$ \A \x^m = \sum_{k=1}^n x_k \A_1 \x^{m-1}=\left(\sum_{k=1}^n x_k\right)^2 \A_{11} \x^{m-2}, $$
where $\A_{11}$ is the root tensor of $\A_1$. The conclusions
(a)-(c) follow immediately. \ep

\subsection{An Algorithm and Numerical Tests}
\hspace{4mm} In Subsection 5.1, we show that a circulant tensor with a diagonal root tensor is positive semi-definite if and only if the global optimal value of (20) or (21) is nonnegative. In this subsection, we present an algorithm to solve the minimization problem (21). Here, $m$ is even, and the norm $\|\cdot\|$ in this section is the $2$-norm.

Suppose $\A \in C_{m,n}$. The minimization problem
\begin{eqnarray}\label{e8.10}
& \min & \A \x^m \nonumber \\
& {\rm subject\ to}\ & \| \x \| = 1,
\end{eqnarray}
can be equivalent to be written as
\begin{eqnarray}\label{e8.11}
& \min & \A \x^1 \cdots \x^m \nonumber \\
& {\rm subject\ to}\ &  \sum_{j=1}^m A_j \x^j =\0  \nonumber   \\
&                    & \| \x^k \| = 1,  \quad  k \in [m],
\end{eqnarray}
where
$$
A_1=\left(\begin{array}{c}
I_m      \\
\vdots   \\
\vdots   \\
-I_m   \end{array}
\right),    \quad
A_2=\left(\begin{array}{c}
-I_m     \\
I_m      \\
\vdots   \\
\vdots   \end{array}
\right),     \quad
\cdots ,     \quad
A_m=\left(\begin{array}{c}
\vdots   \\
\vdots   \\
-I_m     \\
I_m   \end{array}
\right).
$$
Denote $f(\x^1,\cdots,\x^m): =\A \x^1 \cdots \x^m $. Then the augmented Lagrangian function of (\ref{e8.11}) $L_{\beta}(\x^1,\cdots,\x^m,\mathbf{\lambda})$ is defined as
$$
L_{\beta}(\x^1,\cdots,\x^m,\mathbf{\lambda})=f(\x^1,\cdots,\x^m) -
\mathbf{\lambda}^\top \left(\sum_{j=1}^m A_j \x^j\right) +
\frac{\beta}{2} \left\| \sum_{j=1}^m A_j \x^j   \right\|^2 ,
$$
with the given constant $\beta >0$.

We use the alternating direction method of multipliers to solve (\ref{e8.11}).

\begin{Algorithm} Alternating direction method of multipliers for circulant tensors

\noindent
Step 0. Given $\epsilon >0$, $\w ^0=[(\x^1)^0,\cdots,(\x^m)^0, \mathbf{\lambda} ^0] \in \mathbb{R}^{n} \times \cdots \times \mathbb{R} ^{n} \times \mathbb{R} ^{mn}$ and set $k=0$.

\noindent
Step 1. Generate $\w ^{k+1}$ from $\w ^k$, i.e.,  for $j \in [m]$
\begin{equation}\label{e8.12}
(\x^j)^{k+1} = \underset{\x^\top \x =1}{argmin} \ L_{\beta}[(\x^1)^{k+1},\cdots,(\x^{j-1})^{k+1},\x , (\x^{j+1})^{k}, \cdots , (\x^m)^k, \mathbf{\lambda}^k]
\end{equation}
and
$$
\mathbf{\lambda}^{k+1}=\mathbf{\lambda}^k - \beta \sum_{j=1}^m A_j (\x^j)^{k+1} .
$$

\noindent
Step 2. If $\| \w^{k+1}- \w ^k\| < \epsilon $, stop. Otherwise, $k:=k+1$, go to Step 1.

\end{Algorithm}

Note that the subproblem (\ref{e8.12}) is exactly equivalent to a convex quadratic programming on the unit ball, i.e.,
\begin{eqnarray*}
& \min &  \x^\top \x + b^\top \x \nonumber \\
& {\rm subject\ to}\ & \| \x \| = 1,
\end{eqnarray*}
with a given vector $b$. It is well known that it has a closed form solution. So this algorithm is easily implemented.

Under certain condition, the convergence of the algorithm has also been proved, see \cite{HY,JMZ,LWC}. Though the sequence generated from the algorithm may converge to a KKT point, the following numerical results show that the iterative sequence converges to the global minimal solution with a high probability if we choose the initial point randomly. Note that all the diagonal elements of the root tensor are generated randomly in $[-10,10]$.

\begin{Example}\label{exa5}
A circulant tensor $\A \in C_{4,3}$ is generated from a diagonal root tensor with $c_0 = -4.75046$, $c_1= 3.58365$ and $c_2=8.252$.
\end{Example}

\begin{Example}\label{exa6}
A circulant tensor $\A \in C_{4,4}$ is generated from a diagonal root tensor with $c_0 = 3.30134$, $c_1= -9.68746$, $c_2=2.31954$ and $c_3= 7.60276$.
\end{Example}

In the implementation of Algorithm 1, we set the parameters $\beta =
1.2 $ and $\epsilon = 10^{-6}$. And the initial point is generated
randomly. All codes were written by Matlab R2012b and all the
numerical experiments were done on a laptop with Intel Core i5-2430M
CPU 2.4GHz and 1.58GB memory. The numerical results are reported in
Table 1. In the table, $\bar{k}$, $\bar{t}$ and $\bar{\lambda}$
denote the average number of iteration, average time and average
value derived after 100 experiments. ${\lambda}^*$ means the global
minimal solution derived by the polynomial system solver {\bf
Nsolve} available in Mathematica, provided by Wolfram Research Inc.,
Version 8.0, 2010. The frequency of success is also recorded. If $
\| \lambda - \lambda^* \| \leq 10^{-5}$, we say that the algorithm
can find the global minimal solution of (21) successfully.

\begin{table}[!h] \label{tb10}
\caption{Numerical results for Example \ref{exa5} and Example \ref{exa6}}
\begin{center}
\begin{tabular}{cccccc}
\hline\noalign{\smallskip}
                   & $\bar{k}$  & $\bar{t}$  & $\bar{\lambda}$ & $ {\lambda}^* $ & Frequency of success   \\
\noalign{\smallskip}\hline\noalign{\smallskip}
Example \ref{exa5} &62.73       & 0.35607       &-6.39448           & -6.39448        &  100\%            \\
Example \ref{exa6} &92.49       & 0.52795       &-1.79658           & -1.79658        &  100\%            \\
\noalign{\smallskip}\hline
\end{tabular}
\end{center}
\end{table}

From Table 1, we can see that the alternative direction method of
multiplies can be efficient for solving the minimization problem
(21) in some cases. We also test some problems with larger scale.
However, it may be hard to verify the value derived by the algorithm
since the solver {\bf Nsolve} could not work for larger scale
problems.


\end{document}